\documentclass[12pt, letterpaper]{article}
\usepackage[utf8]{inputenc}

\title{Inverting the Sweep Map on $(2n,n)$-Dyck Paths: A Simple Algorithm}
\author{Erin Milne \thanks {emilne@champlain.edu}}
\date{Champlain College}

\begin{document}

\begin{titlepage}
\maketitle
\end{titlepage}

\begin{abstract}
We introduce a simple, rank-based algorithm for inverting the sweep map on $(2n,n)$-Dyck paths.
\end{abstract}

\section{Introduction}

This paper introduces a simple, rank-based algorithm for inverting the sweep map on $(2n,n)$-Dyck paths. Nathan Williams has recently proven in [3] that the sweep map is bijective through analogy with the general modular sweep map, a function that he defines and presents an inversion algorithm for. Adriano Grasia and Gouce Xin present a geometric algorithm for inverting the sweep map on $(m,n)$-Dyck paths when $m$ and $n$ are coprime in [2]; Xin also presents a search algorithm in [4] for finding a sweep map's rank sequence, which is a key step for inversion. Armstrong, Loehr, and Warrington outline a general inversion strategy based on the bounce path in [1]; however, ours is the first paper to present a clear, simple, and efficient inversion algorithm specifically for the $(2n,n)$ case.

\section{The Dyck Path and the Sweep Map}

A \textbf{lattice path} is a path consisting of north and east steps that begins at $(0,0)$ and ends at some point $(a,b)$. A \textbf{$(2n,n)$-Dyck path} is a lattice path that never drops below the line $y=x$ and in which the north steps have length $2n$ and the east steps have length $n$ for some positive integer $n$. Such a path contains $n$ north steps and $2n$ east steps. We assign to each lattice point $(x,y)$ on the path a \textbf{level} equal to $2ny-nx$. Thus a $(2n,n)$-Dyck path begins and ends at level 0. We measure the levels of the north and east steps by the levels of their south and west endpoints, respectively, as Xin does in [4]. Following this west-south convention, we label a south endpoint with an S and a west endpoint with a W. The sweep map sorts these S's and W's in nondecreasing order by level. The sweep map's output is thus a sequence of S's and W's, which we call $\sigma$ according to Xin's convention in [4]. $\sigma$ itself represents the step sequence of another Dyck path; therefore the sweep map sends Dyck paths to Dyck paths. Armstrong, Loehr, and Warrington prove this in [1]. 

It is imperative that we apply the sweep map by scanning from right to left. If we scan from left to right, we do not get a bijection, as the following theorem shows:

\vspace{5mm}

\textbf{Theorem 2.1}: Apply the sweep map to the set of $(2n,n)$-Dyck paths, using the west-south convention. If the sweep map is applied from left to right, it will not yield a bijection. 
	
\textbf{Proof}: It suffices to show that, for a fixed $n$, two different $(2n,n)$-Dyck paths have the same output sequence $\sigma$ when the sweep map is applied from left to right.

	Consider the Dyck path created as follows: First, take a north step from level 0 to level $2n$. Then take two east steps to reach level 0. Then take another north step to reach level $2n$.  Repeat this process until all steps have been used. All south ends have level 0, while all west ends are positive. Thus the sweep map, when applied from left to right, outputs a $\sigma$ sequence of $n$ S's followed by $2n$ W's. 

	Now consider the Dyck path created as follows: Take a north step from level 0 to level $2n$. Then take an east step to level $n$. Then take another north step to level $3n$. Then take two east steps to reach level $n$ again. Then take another north step to level $3n$. Repeat this process until all steps have been used. Scanning from left to right, the first south end has level 0 and the rest have level $n$. The final east step has a west endpoint with level $n$, so this is the last $n$ we encounter as we scan from left to right. All other west ends have level greater than $n$. Thus the sweep map again outputs a $\sigma$ consisting of $n$ S's followed by $2n$ W's. 

\section{Inverting the Sweep Map on $(2n,n)$-Dyck Paths}

The key to inverting the sweep map lies in finding the level of each S and W from the original Dyck path. Thus we must we create what we call $\tau$, the rank sequence consisting of all the levels in order from smallest to largest. It should be obvious to the reader that $\tau$ must contain at least one 0 along with at least one of each multiple of n up to the largest level. We establish bounds for the largest level below. 

\vspace{5mm}

\textbf{Proposition 3.1}: Let $\tau$ be the rank sequence of a $(2n,n)$-Dyck path. Then $2n$$\leq$max$\tau$$\leq$$2n^2$.

\textbf{Proof}: We first consider the upper bound on max$\tau$. The largest value of max$\tau$ is achieved when we take all $n$ north steps consecutively at the beginning of the path. The nth north step then ends at level $2n^2$.

Now we consider the lower bound on max$\tau$. The Dyck path necessarily begins with a north step, which ends at level $2n$. From here, take two east steps to reach level 0. We cannot take any more east steps at this point, since then we would reach a negative level, so we must take another north step to reach level $2n$ again. Then we can take another two east steps to return to level 0. We repeat this process until all steps have been used. Since we never reach a level higher than $2n$, max$\tau$ = $2n$.

\vspace{5mm}

We begin our inversion algorithm by proving a theorem about the levels of the S's.

\vspace{5mm}

\textbf{Theorem 3.2}: Let $\sigma$ be the output sequence of the sweep map applied from right to left on a $(2n,n)$-Dyck path. If $\sigma_{i}$ = W has level $\tau_{i}$ = $kn$ and $\sigma_{i+1}$ = S, then $\tau_{i+1}$ = $kn$ as well. Moreover, if $\sigma_{i+2}$ = S, then $\tau_{i+2}$ = $kn$ also.

\textbf{Proof}: Let $\sigma_{i}$ = W and $\tau_{i}$ = $kn$. Let $\sigma_{i+1}$ = S and assume $\tau_{i+1}$ = $(k+1)n$. This implies that, as we sweep the Dyck path from right to left, we do not encounter any west ends having level $(k+1)n$ before we encounter $\sigma_{i+1}$ = S having level $(k+1)n$. This means that any west end to the right of $\sigma_{i+1}$ = S has level at most $kn$. Assume there is a west end to the right of $\sigma_{i+1}$ = S with level $kn$. The end immediately to the left of this west end must be a south end with level $(k-2)n$. To the left of this south end there can be at most two consecutive west ends (which takes us back to level $kn$) before there must be another south end, necessarily with level $(k-2)n$ again. But to reach a south end with level $(k+1)n$, we must first reach a west end with level at least $(k+3)n$. This west end must ultimately be preceded on the right at some point by a west end with level $(k+2)n$ and another west end with level $(k+1)n$, contradicting our assumption. Since we are sorting in nondecreasing order, $\sigma_{i+1}$ = S must thus have level $\tau_{i+1}$ = $kn$. This same argument also proves that if $\sigma_{i+2}$ = S, then $\tau_{i+2}$ = $kn$ as well.

\vspace{5mm}

We are now ready to present and justify our sweep map inversion algorithm for $(2n,n)$-Dyck paths.

\vspace{5mm}

\textbf{Theorem 3.3}: Inversion Algorithm for the Sweep Map on $(2n,n)$-Dyck Paths:

Input: $\sigma$ (sweep map's output sequence of S's and W's)

Output: A rank sequence $\tau$

Algorithm:

Assign $\tau_{1}=0$

For $i$ = 2 to $2n$:

if $\tau_{i}$ is empty, then:

\hspace{1cm}{ if $\sigma_{i}$ = S, then $\tau_{i}$ = $\tau_{i-1}$ 

\hspace{1cm} else if $\sigma_{i}$ = W, find the number of levels thus far in $\tau$ equal to $\tau_{i-1}$ and then subtract the number of levels thus far in $\tau$ equal to $\tau_{i-1} - 2n$ that correspond to S's. If the difference $x$ is positive, assign level $\tau_{i-1} + n$ to the positions 
$\tau$ corresponding to the positions in $\sigma$ of the following $x$ W's (beginning with $\sigma_{i}$ = W). If $x$ is negative, assign level $\tau_{i-1}$ to the positions in $\tau$ corresponding to the positions in $\sigma$ of the following $x$ W's.

else if $\tau_{i}$ is not empty, return.

\textbf{Proof}: We know that the first level is necessarily a 0, and we have already proven that any south end must have the same level as the end preceding it in $\sigma$. Now consider the W's. If a W has level $kn+n$, the east end of the corresponding east step has level $kn$. This east end could be either an S or another W. Thus we find all elements of $\sigma$ that have level $kn$. Now, each of the elements could have a corresponding west end with level $kn+n$ or a corresponding south end with level $kn-2n$. Thus if we subtract the number of S's with level $kn-2n$ from the total number of ends with level $kn$ and get a positive difference $x$, then we know there are $x$ west ends with level $kn+n$. If $x$ is negative, then the next $x$ W's have level $kn$, since $\tau$ is nondecreasing. Since we can fill multiple positions in $\tau$ at once this way, we make sure that our algorithm only looks at positions in $\tau$ that are currently empty as we loop over $i$ to avoid overwriting positions that have already been filled. 

\vspace{5mm}

We demonstrate our algorithm on a very small example: a (3,6)-Dyck path. Let $\sigma$ = SWSSWWWWW. Set $\tau_{1}$ = 0. $\sigma_{2}$ = W, so we observe that so far there is one 0 and no -6's in $\tau$. Thus we have one W with level 3, so $\tau_{2}$ = 3. $\sigma_{3}$ = S and $\sigma_{4}$ = 3, so $\tau_{3}$ = 3 and $\tau_{4}$ = 3. $\sigma_{5}$ = W, so observe that so far there are three 3's and no -3's in $\tau$. Thus the next three W's have level 6, so $\tau_{5}$ = 6, $\tau_{6}$ =6, and $\tau_{7}$ = 6. $\sigma_{8}$ = W, so we observe that so far there are three 6's in $\tau$ and one 0 that corresponds to an S. Thus the next two W's have level 9, so $\tau_{8}$ = 9 and $\tau_{9}$ = 9. Thus $\tau$ = [0, 3, 3, 3, 6, 6, 6, 9, 9].

We will now use this example to show how to invert the sweep map once we have $\tau$. $\tau_{1}$ corresponds to an S, so we take a north step from level 0 to level 6. We now search for the right-most 6 in $\tau$; since we applied the sweep map by scanning from right to left, the right-most endpoint with a given level corresponds to the left-most endpoint with that level on the original path. The right-most 6 corresponds to a W, so we take an east step to level 3. The right most 3 corresponds to an S, so we take a north step to level 9. The right-most 9 corresponds to a W, so we take an east step to level 6. The second-right-most 6 corresponds to a W, so we take an east step to level 3. The second-right-most 3 corresponds to an S, so we take a north step to level 9. The left-most 9 corresponds to a W, so we take an east step to level 6. The left-most 6 corresponds to corresponds to a W, so we take an east step to level 3. The left-most 3 corresponds to a W, so we take an east step to level 0 and conclude the path. Thus the original Dyck path has step sequence $\sigma^{-1}$ = SWSWWSWWW.

\section{Conclusion}
We are now working on implementing this algorithm in Sage. We are also looking to extend it to $(kn,n)$-Dyck paths, where $k$ is arbitrary. 

\section{Acknowledgements}
I would like to thank Dr. Gregory Warrington of the University of Vermont for introducing me to sweep map research shortly after I earned my master's degree from UVM in May 2015. His unwavering support and guidance have made this work possible. I also thank Dr. Melanie Brown and Dr. Scott Stevens of Champlain College, where I am currently employed as an adjunct math instructor, for their moral support.

\vspace{5mm}

\begin{center}

\textbf{References}

\end{center}

[1] Drew Armstrong, Nicholas A. Loehr, and Gregory S. Warrington, Sweep Maps: A Continuous Family of Sorting Algorithms, \textit{Advances in Mathematics}, 284 (2015) 159-185

[2] Adriano Garsia and Gouce Xin, Inverting the Rational Sweep Map, preprint, arXiv:1602.02346

[3] Nathan Williams, Sweeping Up Zeta, preprint, arXiv:1512.01483

[4] Gouce Xin, An Efficient Search Algorithm for Inverting the Sweep Map on Rational Dyck Paths, preprint, arXiv:1505.00823

\end{document}